\documentclass[12pt]{amsart}
\usepackage{amsmath, amssymb, amsfonts, amscd, amsbsy, amsthm}
\usepackage{latexsym, color, xcolor, graphicx}
\usepackage{bm}
\usepackage{ulem}
\usepackage{mathtools}

\numberwithin{equation}{section}

\theoremstyle{plain}
\newtheorem{Theorem}{Theorem}[section]
\newtheorem{Proposition}[Theorem]{Proposition}

\newtheorem{Problem}[Theorem]{Problem}

\theoremstyle{definition}

\theoremstyle{remark}
\newtheorem{Remark}{{\bf Remark}}

\newcommand{\R}{{\mathbb R}}

\newcommand{\emi}{\textcolor[rgb]{0.00, 0.00, 1.00}}
\newcommand{\emii}{\textcolor[rgb]{1.00, 0.00, 0.00}}

\begin{document}

\title{Maximization of the first Laplace eigenvalue of a finite graph II}

\author[Gomyou]{Takumi Gomyou}
\thanks{Academic Support Center, Kogakuin University, 
Nakano-machi, Hachioji 192-0015, Japan, kt13804@ns.kogakuin.ac.jp}

\author[Nayatani]{Shin Nayatani}
\thanks{Graduate School of Mathematics, Nagoya University,
Chikusa-ku, Nagoya 464-8602, Japan, nayatani@math.nagoya-u.ac.jp}



\maketitle

\begin{abstract}
Given a length function on the set of edges of a finite graph, the corresponding Fujiwara Laplacian is defined.
We consider a problem of maximizing the first nonzero eigenvalue of this graph Laplacian over all choices of edge-length function subject to a certain normalization.
In this paper we prove that the supremum of the first nonzero eigenvalue is finite if and only if the graph is a tree. 
We also prove that the supremum of the first nonzero eigenvalue is nonincreasing under taking a subgraph. 
\end{abstract}

\section*{Introduction}

This paper is a sequel to \cite{GomyouNayatani}, where we initiated the study of a maximization problem for  the first nonzero eigenvalue of a certain Laplacian on a finite graph.

Let $G=(V,E)$ be a finite graph, where $V$ (resp. $E$) is the set of vertices (resp. edges).
If an edge-length function $l\colon E\to \R_{>0}$ is given, one can define the corresponding {\it Fujiwara Laplacian} $\Delta_l$ \cite{Fujiwara}; This is the Laplacian of our interest.
Then we regard the length function $l$ as a variable and consider the problem of maximizing the first nonzero eigenvalue of the Laplacian $\Delta_l$ over all choices of $l$ subject to a certain normalization.
The main result of \cite{GomyouNayatani} states that if a maximizing (or more generally, extremal) length function exists, then there exists a map $\varphi\colon V \to \R^N$ consisting of first eigenfunctions of the corresponding Laplacian so that the length function can be explicitly expressed in terms of the map $\varphi$ and the Euclidean distance.
Here, $N$ is some positive integer less than or equal to the multiplicity of the first nonzero eigenvalue.
In \cite{GomyouNayatani}, we solved the maximization problems for all graphs with up to four vertices by numerical means.
In particular, we determined the supremum of the first nonzero eigenvalue, written as $\Lambda_1(G)$, for these graphs and found that $\Lambda_1(G)$ is finite if and only $G$ is a tree, that is, $G$ contains no nontrivial cycles. 

One of the main results of this paper asserts that this equivalence holds true completely in general.
Namely, we prove
\begin{Theorem}\label{main1}
Let $G$ be a finite connected graph. Then $\Lambda_1(G)$ is finite if and only if $G$ is a tree. 
\end{Theorem}

We prove the only-if part of Theorem \ref{main1} by showing its contraposition, that is, that if $G$ contains 
a nontrivial cycle, then $\Lambda_1(G)=\infty$. 
In fact, if we choose the length of one edge of any cycle almost equal to $1/2$ and let the lengths 
of other edges tend to zero at the same speed, then we can verify that the first eigenvalue diverges to infinity. 
The if part of Theorem \ref{main1} is proved by evaluating the Rayleigh quotient for a test function with suitable properties. 
In fact, this gives an explicit upper bound of $\Lambda_1(G)$ in terms of the number of vertices. 
In the particular cases of the path graphs $P_n$ on $n$ vertices, we observe that a natural choice of test function provides a uniform bound which does not depend on $n$. 

The other of the main results establishes the monotonicity of the invariant $\Lambda_1$. Namely, we prove

\begin{Theorem}\label{main2}
If $G$ is a finite connected graph and $H$ is a subgraph of $G$, then the inequality 
$\Lambda_1(H)\leq \Lambda_1(G)$ holds.\end{Theorem}

In order to prove Theorem \ref{main2}, we first prove that if the graph $G$ has a pendant edge, that is, an edge incident to a degree-one vertex, and $G'$ is the graph obtained by contracting the pendant edge, then they satisfy $\Lambda_1(G')\leq \Lambda_1(G)$. If $G$ has no pendant edge, we cut $G$ at a vertex and produce $G_{\mathrm{cut}}^{\mbox{}}$ with a pendant edge, and observe that $\Lambda_1(G_{\mathrm{cut}}^{\mbox{}})\leq \Lambda_1(G)$ holds. 
We arrive at the conclusion by repeating these procedures suitably until $G$ is collapsed to $H$.

The problem of maximizing the first Laplace eigenvalue on a graph was first considered by Fiedler \cite{Fiedler2}. More recently, Fiedler's problem was generalized by G\"{o}ring-Helmberg-Wappler \cite{GoringHelmbergWappler1, GoringHelmbergWappler2}. They fix a vertex-weight $m_0$ and an edge-length function $l$, and maximize the first nonzero eigenvalue of the weighted Laplacian $\Delta_{(m_0,m_1)}$ over all choices of edge-weight $m_1$ subject to a normalization which concerns $l$.
It is worth mentioning that in their eigenvalue maximization problem the supremum of the first nonzero eigenvalue in question is necessarily finite.

In the background of the problem considered in this paper, there also is the problem of maximizing the first Laplace eigenvalue on a manifold. In fact, Berger \cite{Berger} asked whether on an arbitrary compact manifold of dimension $n$, the scale-invariant quantity $\lambda_1(g) \mathrm{Vol}(g)^{2/n}$, where $g$ is a Riemannian metric and $\mathrm{Vol}(g)$ is the volume of $g$, was bounded from above by a constant depending only on the manifold.
Urakawa \cite{Urakawa} then answers Berger's question in the negative by giving an explicit family of metrics $g_t$, $t>0$, on the three-sphere for which the quantity $\lambda_1(g_t) \mathrm{Vol}(g_t)^{2/3}$ diverges to infinity as $t\to\infty$. As a corollary of this result and its analogues for higher dimensional spheres, it is now known that the quantity $\lambda_1(g) \mathrm{Vol}(g)^{2/n}$ is unbounded on any compact manifold of dimension $\geq 3$.
On the other hand, for surfaces, the quantity $\lambda_1(g) \mathrm{Area}(g)$ is bounded from above and rich progress has been made on the problem of maximizing this quantity.
See \cite{LiYau, Nadirashvili, Petrides, MatthiesenSiffert, Ros, Karpukhin}.

This paper is organized as follows.
In Section 1, we review the definition of the Fujiwara Laplacian and the formulation \cite{GomyouNayatani} of the maximization problem for the first nonzero eigenvalue of this Laplacian.
In Section 2, we prove that if a finite connected graph contains a nontrivial cycle, then the supremum of the first eigenvalue under the length normalization diverges.
In Section 3, we show that the supremum of the first eigenvalues for a tree is bounded above by a constant depending on the number of vertices, thereby concluding the proof of the first main theorem. 
In Section 4, we obtain a sharper bound, which does not depend on the number of vertices, for path graphs.
In Section 5, we investigate the limiting behavior of eigenvalues of the Laplacian when a pendant edge 
is contracted . We then prove the second main theorem. 

\section{Preliminaries}
Let $G=(V, E)$ be a finite connected graph, where $V$ and $E$ are the sets of vertices and (undirected) edges, respectively. We assume that $G$ is simple, that is, that $G$ has no loops nor multiple edges. Let $uv$ denote the edge whose endpoints are $u$ and $v$. $G$ being undirected means that $uv=vu$. 
If a vertex weight $m_0\colon V \to \R_{>0}$ and an edge weight $m_1\colon E \to \R_{\geq 0}$ are chosen, 
the corresponding graph Laplacian $\Delta_{(m_0,m_1)}\colon \R^V \to \R^V$ is defined by 
\begin{equation*}
(\Delta_{(m_0,m_1)} \varphi) (u) = \sum_{v\sim u} \frac{m_1(uv) }{m_0(u)} \left( \varphi(u) - \varphi(v) \right),\quad u\in V,
\end{equation*}
where we write $v\sim u$ if $uv\in E$. 
Here, $\R^V$ denote the set of functions $\varphi\colon V\to \R$, equipped with the inner product
$$
\langle \varphi_1, \varphi_2 \rangle = \sum_{u\in V} m_0(u)\, \varphi_1(u)\, \varphi_2(u),
\quad \varphi_1, \varphi_2\in \R^V.
$$
Note that $\Delta_{(m_0,m_1)}$ is a nonnegative symmetric linear operator on $\R^V$. 

Given an edge-length function $l\colon E \to \R_{>0}$, Fujiwara \cite{Fujiwara} defined a vertex-weight $m_0$ and an edge-weight $m_1$ by
$$
m_0(u) = \sum_{v \sim u} l(uv),\,\, u \in V\quad \mbox{and}\quad
m_1(uv) = \frac{1}{l(uv)},\,\, uv \in E,
$$
respectively.
He then considered the corresponding Laplacian.
We denote this Laplacian by $\Delta_l$ and call it the {\it Fujiwara Laplacian}.
Note that $\Delta_l\colon \R^V \to \R^V$ is a positive symmetric linear operator, hence has only nonnegative real eigenvalues; It always has eigenvalue $0$, and the corresponding eigenspace consists precisely of constant functions on $V$.
The second smallest eigenvalue of $\Delta_l$, which is positive, will be denoted by $\lambda_1(G, l)$
and referred to as the first nonzero eigenvalue of $\Delta_l$.
It is a standard fact that $\lambda_1(G, l)$ is characterized variationally as
\begin{eqnarray}\label{lambda1-characterization1}
\lambda_1(G, l) &=&
\min_{\varphi} \frac{\langle \Delta_l \varphi, \varphi \rangle}{
\langle \varphi, \varphi \rangle} \nonumber\\
&=&
\min_{\varphi} \frac{\sum_{uv\in E} m_1(uv)
( \varphi(u) - \varphi(v) )^2}{\sum_{u\in V} m_0(u) 
\varphi(u)
^2},
\end{eqnarray}
where 
the minimum is taken over all nonzero 
functions $\varphi$ such that\\ $\sum_{u\in V} m_0(u) \varphi(u)=0$, meaning that
$\varphi$ is orthogonal to constant functions, that is, eigenfunctions of the eigenvalue $0$.

Other than the operator language, we can employ the matrix language to describe the graph Laplacian.
To do so, we write $V=\{ 1,\dots, |V| \}$ and choose the orthonormal basis
$$
e_i\colon j\in V \mapsto \frac{\delta_{ij}}{\sqrt{m_0(i)}}\in \R,\quad i\in V,
$$
of $\R^V$.
Then the corresponding representation matrix $L := L_l$ of the Laplacian $\Delta_l$
is given by $L= D^{-1/2} L_0 D^{-1/2}$, where $D=\mathrm{diag}(m_0(1),\dots, m_0(|V|))$ and $L_0$
has diagonal components
$$
(L_0)_{i,i}=\sum_{j\sim i} m_1(ij),\quad i\in V
$$
and off-diagonal components
$$
(L_0)_{i,j}= \left\{
\renewcommand{\arraystretch}{1.5}
\begin{array}{cc} - m_1(ij),
& 
ij\in E,\\
0, & ij\notin E.
\end{array}
\right.
$$
The matrix $L$ is positive, symmetric, and has eigenvalue $0$ with eigenvector
$\bm{x}_0 = {}^t \bigl(\sqrt{m_0(1)},\dots,\sqrt{m_0(|V|)}\bigr)$.
Note also that the variational characterization \eqref{lambda1-characterization1} of
$\lambda_1(G, l)$ is expressed as
\begin{equation*}
\lambda_1(G, l) = \min_{\bm{x}} \frac{\langle L\bm{x}, \bm{x} \rangle}{\langle \bm{x}, \bm{x} \rangle},
\end{equation*}
where $\langle \cdot, \cdot \rangle$ denotes the Euclidean inner product on $\R^{|V|}$ and the minimum is taken over
all nonzero vectors $\bm{x} \in \R^{|V|}$ such that $\langle \bm{x}, \bm{x}_0 \rangle = 0$.

We now recall the first-eigenvalue maximization problem \cite{GomyouNayatani} whose variable is the edge-length function $l$.

\begin{Problem}\label{maxspecGN}
Over all edge-length functions $l$, subject to the normalization
\begin{equation}\label{constraint}
\sum_{u \in V} m_0(u) = 2 \sum_{uv\in E} l(uv) = 1,
\end{equation}
maximize the first nonzero eigenvalue $\lambda_1(G,l)$ of $\Delta_{l}$.
\end{Problem}

Set $\Lambda_1(G) := \sup_l \lambda_1(G,l)$, where $l$ runs over all edge-length functions
satisfying \eqref{constraint}.
Note that it is also expressed as
\begin{equation}\label{Lambda1}
\Lambda_1(G) = \sup_l \lambda_1(G,l)\cdot \Bigl( \sum_{u \in V} m_0(u) \Bigr)^2
\end{equation}
with $l$ unnormalized.

\section{Divergence for a graph containing a cycle}

In this section, we prove the only-if part of Theorem \ref{main1} by proving its contraposition. 

\begin{Theorem}
Let $G$ be a finite connected graph which contains a nontrivial cycle.
Then $\Lambda_1(G)=\infty$.
\end{Theorem}

\begin{proof}
Suppose that $G$ contains a nontrivial cycle $C$, and choose an edge $e_0$ of $C$.
The subgraph $G \setminus \{ e_0 \}$ remains connected. 
By renumbering of vertices, we suppose that 
$e_0=(n-1, n)$. 
For small $t>0$, define a length function $l_t\colon E\to \R_{>0}$ by
$$
l_t(e) = \left\{\begin{array}{cl}
1, & e = e_0, \\
t, & e \neq e_0.
\end{array}
\right.
$$
Therefore, the corresponding edge-weight $m_1\colon E\to \R_{\geq 0}$ is 
$$
m_1(e) = \left\{\begin{array}{cl}
1, & e = e_0, \\
1/t, & e \neq e_0.
\end{array}
\right.
$$
Note that the vertex-weight $m_0\colon V\to \R_{>0}$ is given by
$$
m_0(i) = \left\{ \begin{array}{cl} \deg(i)\,t, &1\leq i\leq n-2,\\ 1+ (\deg(i)-1)\,t, & i=n-1, n, \end{array} \right.
$$
and satisfies  
$$
\sum_{i=1}^n m_0(i) = 2(1+(|E|-1)t) = 2+O(t).
$$
It is therefore sufficient to verify that $\lambda_1(G,l_t)$ diverges to positive infinity as $t\to 0$ 
in order to prove the theorem. 

The unnormalized Laplacian matrix of $G$ with respect to the edge-weight $m_1$ is given by $$L_0(G,m_1) = \frac{1}{t} L_0(G \setminus \{ e_0 \}, {\bm 1}) + 
\left(\begin{array}{c|c}
\mbox{\huge 0} & \mbox{\huge 0} \\ \hline
\raisebox{-0.5em}{\mbox{\huge 0}} &
\begin{array}{cc}
1 & -1 \\
-1 & 1
\end{array}
\end{array}\right).$$
Let $\bm x=(x_1,\dots,x_n)$ be an eigenvector belonging to the first nonzero eigenvalue of $L_0(G,m_1)$. Note that ${\bm x}$ is perpendicular to the vector $(1,\dots,1)$. Then 
\begin{eqnarray*}
{}^t {\bm x} L_0(G,m_1) {\bm x} &=& \frac{1}{t} {}^t {\bm x} L_0(G \setminus \{ e_0 \},{\bm 1}) {\bm x} 
+ (x_{n-1}-x_n)^2\\ 
&\geq& \frac{1}{t} {}^t {\bm x} L_0(G \setminus \{ e_0 \}, {\bm 1}){\bm x}. 
\end{eqnarray*}
Therefore, 
$$
\lambda_1(L_0(G,m_1)) = \frac{{}^t {\bm x} L_0(G,m_1) {\bm x}}{{}^t {\bm x}{\bm x}} 
\geq \frac{1}{t} \frac{{}^t {\bm x} L_0(G \setminus \{ e_0 \}, {\bm 1}){\bm x}}{{}^t {\bm x}{\bm x}} 
\geq \frac{1}{t} \lambda_1(L_0(G \setminus \{ e_0 \}, {\bm 1})).
$$
The connectedness of $G \setminus \{ e_0 \}$ assures $\lambda_1(L_0(G \setminus \{ e_0 \}, {\bm 1}))>0$, and therefore, the right-hand side diverges to positive infinity as $t\to 0$. 
Since $L(G,m_1) = D^{-1/2} L_0(G,m_1) D^{-1/2}$ and $D^{-1/2}\geq c I$, where $c$ is a positive constant, $\lambda_1(G,l_t)=\lambda_1(L(G,m_1))$ also diverges to positive infinity as $t\to 0$. 

\end{proof}

\section{Finiteness for a tree}
In this section, we complete the proof of Theorem \ref{main1} by showing that $\Lambda_1(G)$ 
is finite if $G$ is a tree. 

\begin{Theorem}
Let $T$ be a finite connected tree with $n$ vertices. 
Then $\Lambda_1(T)\leq 16(n-1)^3$.
\end{Theorem}

\begin{proof}
Let $T = (V,E)$ be a tree and write $V=\{1,\dots,n\}$.
Let $l$ be any edge-length function satisfying the normalization condition $\sum_{ij \in E} l(ij)=1/2$.
Let $i_0j_0$ be an edge of maximal length, i.e., $l(i_0j_0) = \max_{ij \in E} l(ij)$. 
It is easy to see that there exists a map $\varphi\colon V\to \R$ such that $|\varphi(i)-\varphi(j)|=l(ij)$ for each edge $ij\in E$ and $\sum_{i=1}^nm_0(i)\varphi(i) = 0$.
Then we have
$$
\lambda_1(T, l)\leq \frac{\sum_{ij \in E} m_1(ij)(\varphi(i)-\varphi(j))^2}{\sum_{i=1}^n
m_0(i) \varphi(i)^2}.
$$
Since $m_1(ij)=1/l(ij)$, the numerator of the right-hand side is 
$$
\sum_{ij \in E} \frac{1}{l(ij)} (\varphi(i)-\varphi(j))^2 = \sum_{ij \in E} l(ij) = \frac{1}{2}. 
$$
We estimate the denominator as 
\begin{eqnarray*}
\sum_{i=1}^n m_0(i) \varphi(i)^2 &\geq& m_0(i_0) \varphi(i_0)^2 + m_0(j_0) \varphi(j_0)^2\\
&\geq& l(i_0j_0)\cdot \left( \frac{l(i_0j_0)}{2} \right)^2\,\, =\,\, \frac{l(i_0j_0)^3}{4}\\
&\geq& \frac{1}{4} \left( \frac{1}{2(n-1)} \right)^3\,\, =\,\, \frac{1}{32(n-1)^3}.
\end{eqnarray*}
Here, we have used $\max \{ | \varphi(i_0) |, | \varphi(j_0) | \}\geq \frac{l(i_0j_0)}{2}$ and
$l(i_0j_0)\geq \frac{1}{2|E|} = \frac{1}{2(n-1)}$.
Therefore, $\lambda_1(T, l)\leq 16(n-1)^3$. 
\end{proof}

\section{Sharper upper bound for path graphs}
As shown in the previous section, the invariant $\Lambda_1$ for a tree is bounded above by a cubic polynomial in the number of vertices. 
It is manifest from the proof that this upper bound should be far from being tight. 
In fact, for a path graph, we can obtain a much sharper bound.
\begin{Theorem}
Let $P_n$ be the path graph with n vertices.
Then $\Lambda_1(P_n) \leq 24$.
\end{Theorem}
\begin{proof}
Let $P_n$ be the path graph on $n$ vertices $1, \cdots, n$, where the vertex $i$ is adjacent to $i+1$ for $1\leq i\leq n-1$. 
Choose any edge-length function and denote the length of the edge $(j,j+1)$ by $l_j$ for $1\leq j\leq n-1$. 
Define a map $\varphi$ by $\varphi(1)=0$ and $\varphi(i)=\sum_{j=1}^{i-1} l_j$.
We consider a problem of maximizing
\begin{equation}\label{RayleighQuotient}
\frac{\sum_{j=1}^{n-1} m_1(j, j+1) (\varphi(j)-\varphi(j+1))^2}{\sum_{i=1}^n m_0(i) (\varphi(i)-\overline{\varphi})^2}
\end{equation}
over all edge-length functions $l=(l_j)$ subject to $\sum_{j=1}^{n-1} l_j = 1/2$ and $l_j > 0$, $1 \leq j \leq n-1$.
The maximum value of this problem provides an upper bound for $\Lambda_1(P_n)$.
It can be verified that the barycenter $\overline{\varphi}$ is constant and equal to $1/4$. In fact, 
\begin{eqnarray*}
\overline{\varphi} &=& \sum_{i=1}^{n}m_0(i)\varphi(i)\\
&=& \sum_{i=2}^{n-1} (l_{i-1} +l_i) \sum_{j=1}^{i-1} l_j + l_{n-1} \sum_{j=1}^{n-1} l_j \\
&=& \sum_{i=1}^{n-2} l_i \sum_{j=1}^i l_j +\sum_{i=2}^{n-1} l_i \sum_{j=1}^{i-1} l_j+ l_{n-1} \sum_{j=1}^{n-1} l_j \\
&=& \left( \sum_{i=1}^{n-2} l_i^2 +\sum_{i=2}^{n-2} l_i \sum_{j=1}^{i-1} l_j \right) +\sum_{i=2}^{n-1} l_i \sum_{j=1}^{i-1} l_j +\left(l_{n-1} \sum_{j=1}^{n-2} l_j +l_{n-1}^2 \right) \\
&=& \sum_{i=1}^{n-1} l_i^2 +2\sum_{i=2}^{n-1} l_i \sum_{j=1}^{i-1} l_j \\ 
&=& \left( \sum_{i=1}^{n-1} l_i \right)^2\\ 
&=& \frac{1}{4}.
\end{eqnarray*}
Since the numerator of \eqref{RayleighQuotient} is constant and equal to $1/2$ and the constraint allows to eliminate the variable $l_{n-1}$ by $l_{n-1} = 1/2 -\sum_{j=1}^{n-2} l_j$, the problem reduces to the minimization of $F(l) = \sum_{i=1}^{n} m_0(i) (\varphi(i)-\overline\varphi)^2$ 
over the range $\mathcal{D} =\{ l = (l_j)_{1 \leq j \leq n-2} \mid l_j > 0, 1 \leq j \leq n-2 \ \text{and} \ \sum_{j=1}^{n-2} l_j < 1/2 \}$. 
The objective function $F$ becomes 
\begin{eqnarray*}
F(l) 
&=& \sum_{i=1}^{n}m_0(i)\varphi(i)^2 -\sum_{i=1}^{n} m_0(i)\overline\varphi^2 \\
&=& \sum_{i=1}^{n-2}m_0(i)\varphi(i)^2 +(l_{n-2}+l_{n-1})\left( \sum_{j=1}^{n-2}l_j \right)^2 +l_{n-1}\left( \sum_{j=1}^{n-1}l_j \right)^2 -\frac{1}{16} \\
&=& \sum_{i=1}^{n-2}m_0(i)\varphi(i)^2 +\left( \frac{1}{2} -\sum_{i=1}^{n-3}l_i \right) \left( \sum_{j=1}^{n-2}l_j \right)^2 +\left( \frac{1}{2} -\sum_{j=1}^{n-2}l_j \right) \frac{1}{4} -\frac{1}{16}.
\end{eqnarray*}
The partial derivatives of $F$ with respect to $l_k$, $1\leq k\leq n-3$, are  
\begin{eqnarray*}
\frac{\partial F}{\partial l_k}
&=& \varphi(k)^2 +\varphi(k+1)^2 +2\sum_{i=k+1}^{n-2} m_0(i)\varphi(i) -\left( \sum_{j=1}^{n-2} l_j \right)^2 +2\left( \frac{1}{2} -\sum_{i=1}^{n-3}l_i \right) \sum_{j=1}^{n-2} l_j -\frac{1}{4}.
\end{eqnarray*}
The sum of the first three terms is
\begin{eqnarray*}
&& \left(\sum_{j=1}^{k-1} l_j \right)^2 +\left(\sum_{j=1}^{k-1} l_j +l_k \right)^2 
+2\sum_{i=k+1}^{n-2} (l_{i-1} +l_i) \sum_{j=1}^{i-1} l_j \\
&=& l_k^2 +2\left(\sum_{j=1}^{k-1} l_j \right)^2 +2l_k\sum_{j=1}^{k-1} l_j  
+2\left\{ \sum_{i=k}^{n-3} l_i \left(\sum_{j=1}^{i-1} l_j +l_i \right) +\sum_{i=k+1}^{n-2} l_i \sum_{j=1}^{i-1} l_j \right\} \\
&=& l_k^2 +2\left(\sum_{j=1}^{k-1} l_j \right)^2 +2l_k\sum_{j=1}^{k-1} l_j 
+2\left\{ l_k \sum_{j=1}^{k-1} l_j +2\sum_{i=k+1}^{n-3} l_i \sum_{j=1}^{i-1} l_j +\sum_{i=k}^{n-3} l_i^2 +l_{n-2} \sum_{j=1}^{n-3} l_j \right\} \\
&=& l_k^2 +2\left(\sum_{j=1}^{k-1} l_j \right)^2 +4l_k\sum_{j=1}^{k-1} l_j 
+2\left\{ 2\sum_{i=k+1}^{n-3} l_i \sum_{j=1}^{k-1} l_j +2\sum_{i=k+1}^{n-3} l_i \sum_{j=k}^{i-1} l_j +\sum_{i=k}^{n-3} l_i^2 +l_{n-2} \sum_{j=1}^{n-3} l_j \right\} \\
&=& l_k^2 +2\left(\sum_{j=1}^{k-1} l_j \right)^2 +4\sum_{i=k}^{n-3} l_i \sum_{j=1}^{k-1} l_j +2\left( \sum_{i=k}^{n-3} l_i \right)^2 +2l_{n-2} \sum_{j=1}^{n-3} l_j \\
&=& l_k^2 +2\left( \sum_{j=1}^{n-3} l_j \right)^2 +2l_{n-2} \sum_{j=1}^{n-3} l_j.
\end{eqnarray*}
Therefore, we conclude 
\begin{eqnarray*}
\frac{\partial F}{\partial l_k} 
&=& l_k^2 +2\left( \sum_{j=1}^{n-3} l_j \right)^2 +2l_{n-2} \sum_{j=1}^{n-3} l_j -\left( \sum_{j=1}^{n-2} l_j \right)^2 +\sum_{j=1}^{n-2} l_j -2\sum_{i=1}^{n-3}l_i \left( \sum_{j=1}^{n-3} l_j +l_{n-2} \right) -\frac{1}{4} \\
&=& l_k^2 -\left( \sum_{j=1}^{n-2} l_j \right)^2 +\sum_{j=1}^{n-2} l_j -\frac{1}{4}.
\end{eqnarray*}
The partial derivative of $F$ with respect to $l_{n-2}$ is 
\begin{eqnarray*}    
\frac{\partial F}{\partial l_{n-2}} &=& \varphi(n-2)^2 + 2 \left(\frac{1}{2} - \sum_{i=1}^{n-3} l_i \right) \sum_{j=1}^{n-2} l_j -\frac{1}{4} \\
&=& -\left(\sum_{j=1}^{n-3} l_j \right)^2 +\sum_{j=1}^{n-2} l_j -2 l_{n-2} \sum_{j=1}^{n-3} l_j -\frac{1}{4} \\
&=& l_{n-2}^2 -\left( \sum_{j=1}^{n-2} l_j \right)^2 +\sum_{j=1}^{n-2} l_j -\frac{1}{4}.
\end{eqnarray*}
It follows that the only solution of $\frac{\partial F}{\partial l_k} = 0$ in the range 
$\mathcal{D}$ is $l_k = \frac{1}{2(n-1)}$, $1\leq k\leq n-2$. 
The second-order partial derivatives are
\begin{eqnarray*}
\frac{\partial^2 L}{\partial l_{k'} \partial l_k} &=& 2\delta_{k, k'} \, l_k -2\sum_{j=1}^{n-2}l_j +1, \quad 1 \leq k,k' \leq n-2. 
\end{eqnarray*}
Thus, the Hessian matrix of $F$ is 
\begin{eqnarray*}    
2
\begin{pmatrix}
l_1 & & \text{\huge{0}} \\
& \ddots & \\
\text{\huge{0}} & & l_{n-2}
\end{pmatrix}
+\left( 1 -2\sum_{j=1}^{n-2} l_j \right)
\begin{pmatrix}
1 & \cdots & 1 \\
\vdots & \ddots & \vdots \\
1 & \cdots & 1
\end{pmatrix},
\end{eqnarray*}
and it is positive definite in the range $\mathcal{D}$.
Hence, $F$ is convex in $\mathcal{D}$ and the uniform edge-length function $l_k = \frac{1}{2(n-1)}$ is a unique minimizer of $F$ in the range $\mathcal{D}$. 
For this $l$, the objective value is 
$$\frac{\sum_{ij \in E} (1/l(ij)) (\varphi(i)-\varphi(j))^2}{\sum_{i=1}^n
m_0(i) \varphi(i)^2 -\overline{\varphi}^2} = \frac{24(n-1)^2}{n^2-2n+3}\leq 24.$$ 
\end{proof}

\section{Limit behavior of spectrum}

Let $G$ be a finite connected graph equipped with an arbitrary edge-length function $l$.
Let $u_1,\dots, u_n$ denote the vertices of $G$. 
We consider a graph $\widetilde{G}$ by introducing a new vertex $u_{n+1}$ and joining it to $u_n$.
For $t>0$, we extend $l$ to the edge-length function $l_t$ of $\widetilde{G}$ by setting
the length of the edge $u_n u_{n+1}$ to $t$.

In this section, we investigate the limiting behavior of the Laplace eigenvalues of
$(\widetilde{G}, l_t)$ as $t$ approaches $0$.
If $L=(L_{i,j})$ denotes the Laplacian matrix of $(G, l)$, then that of $(\widetilde{G}, l_t)$ is
given by
$$
\widetilde{L} = \left( \begin{array}{ccccc}
L_{1,1} & \cdots & L_{1,n-1} & L_{1,n}+O(t) & 0\\
\vdots & \ddots &\vdots & \vdots & \vdots \\
L_{n-1,1} & \cdots & L_{n-1,n-1} & L_{n-1,n}+O(t) & 0\\
L_{n,1}+O(t) & \cdots & L_{n,n-1}+O(t) & \widetilde{L}_{n,n} & \widetilde{L}_{n,n+1}\\
0 &\cdots & 0 & \widetilde{L}_{n+1,n} & 1/t^2
\end{array} \right),
$$ 
where
$$
\widetilde{L}_{n,n} = L_{n,n} + \frac{1}{m_0(n)} \frac{1}{t} \left( 1 - \frac{t}{m_0(n)} + O(t^2) \right)
$$
and
$$
\widetilde{L}_{n,n+1} = \widetilde{L}_{n+1,n}
= - \frac{1}{\sqrt{m_0(n)}} \frac{1}{t^{3/2}} \left( 1 -\frac{t}{2m_0(n)} +O(t^2) \right).
$$

\begin{Theorem}\label{contract}
As $t$ approaches $0$:
\begin{enumerate}
\renewcommand{\theenumi}{\roman{enumi}}
\renewcommand{\labelenumi}{(\theenumi)}
\item The largest eigenvalue of $\widetilde{L}$ diverges to positive infinity.
\item All the other eigenvalues of $\widetilde{L}$ converge to those of $L$.
\end{enumerate}
In particular, we have
$$
\lim_{t\to 0} \lambda_1(\widetilde{G}, l_t) = \lambda_1(G, l),
$$
and therefore,
$$
\Lambda_1(\widetilde{G})\geq \Lambda_1(G).
$$
\end{Theorem}

In order to prove the theorem, we consider more general matrices modelled on $L$ and $\widetilde{L}$.
Let $A=(a_{i,j})$ be an $n \times n$ symmetric matrix.
We consider the following $(n+1)\times(n+1)$ matrix:
$$
\widetilde{A}=
\left(
\begin{array}{ccccc}
a_{1,1} & \cdots & a_{1,n-1} & a_{1,n}+O(t) & 0\\
\vdots & \ddots &\vdots & \vdots & \vdots \\
a_{n-1,1} & \cdots & a_{n-1,n-1} & a_{n-1,n}+O(t)& 0\\
a_{n,1}+O(t) & \cdots & a_{n,n-1}+O(t) & a_{n,n} +\frac{t^{-1}}{\alpha}(1-\frac{t}{\alpha}+O(t^2))
& -\frac{t^{-3/2}}{\sqrt{\alpha}}(1-\frac{t}{2\alpha}+O(t^2))\\
0 &\cdots & 0 & -\frac{t^{-3/2}}{\sqrt{\alpha}}(1-\frac{t}{2\alpha}+O(t^2)) & t^{-2}
\end{array}
\right),
$$
where $\alpha$ is a positive constant.

\begin{Proposition}\label{asymptotic}
As $t$ approaches $0${\rm :}
\begin{enumerate}
\renewcommand{\theenumi}{\roman{enumi}}
\renewcommand{\labelenumi}{(\theenumi)}
\item The largest eigenvalue of $\widetilde{A}$ diverges to positive infinity.
\item All the other eigenvalues of $\widetilde{A}$ converge to those of $A$.
\end{enumerate}
\end{Proposition}

\begin{proof}
We prove (i).
Let $\bm{u} = {}^t (0,\cdots, 0, -\sqrt{t/\alpha}, 1)\in\R^{n+1}$.
Then $\langle \bm{u}, \bm{u}\rangle = 1 + O(t)$ and
$\langle\widetilde{A}\bm{u}, \bm{u}\rangle = O(t^{-2})$ as $t\to 0$.
Therefore, the conclusion follows.

Next we prove (ii).
For any $\bm{v} = {}^t (v_1,\cdots, v_n)\in\R^n$, let $\widetilde{\bm{v}} = {}^t (v_1,\cdots, v_{n-1}, v_n,\sqrt{t/\alpha}\, v_n)\in\R^{n+1}$.
Note that $\langle \widetilde{\bm{v}}, \widetilde{\bm{v}} \rangle =
( 1 + O(t)) \langle \bm{v}, \bm{v} \rangle$.
We compute
\begin{eqnarray*}
\widetilde{A}\widetilde{\bm{v}} &=&
\biggl(
\sum_{j=1}^n a_{1,j} v_j + O(t) v_n, \cdots, \sum_{j=1}^n a_{n-1,j} v_j + O(t) v_n,
\sum_{j=1}^n a_{n,j}v_j - \frac{v_n}{2\alpha^2} + O(t) \Bigl(\sum_{i=1}^n v_i\Bigr), \\
&& \frac{v_n}{2\alpha^{3/2}} t^{-1/2} + O(t^{1/2}) v_n \biggr).
\end{eqnarray*}
It follows that
$$
\langle \widetilde{A}\widetilde{\bm{v}}, \widetilde{\bm{v}} \rangle
= \langle A\bm{v}, \bm{v} \rangle +  O(t) \langle \bm{v}, \bm{v} \rangle
\quad \mbox{and}\quad
\frac{\langle \widetilde{A}\widetilde{\bm{v}}, \widetilde{\bm{v}} \rangle}{
\langle \widetilde{\bm{v}}, \widetilde{\bm{v}} \rangle}
= \frac{\langle A\bm{v}, \bm{v} \rangle}{\langle \bm{v}, \bm{v} \rangle} + O(t) .
$$

Let
$$
0=\lambda_0 < \lambda_1 \leq \cdots \leq \lambda_{n-1}\quad \mbox{and}\quad
0=\widetilde{\lambda}_0 < \widetilde{\lambda}_1 \leq \cdots \leq \widetilde{\lambda}_{n-1}
\leq \widetilde{\lambda}_n
$$
denote the eigenvalues of $A$ and $\widetilde{A}$, respectively.
Let $\bm{e}_0, \cdots, \bm{e}_{n-1}$ be orthonormal eigenvectors of $A$ such that
$A \bm{e}_k = \lambda_k \bm{e}_k$ for $0 \leq k \leq n-1$.
Set $V_k := \text{span}\{ \bm{e}_0,\cdots,\bm{e}_{k-1} \}$, so that its orthogonal complement in $\R^n$ is
$V_k^{\perp} =\text{span}\{ \bm{e}_k,\cdots,\bm{e}_{n-1} \}$.
With $\bm{u}$ as above, set
$\widetilde{V}_k :=\text{span}\{ \widetilde{\bm{e}}_k,\cdots,\widetilde{\bm{e}}_{n-1},\bm{u} \}^{\perp}$, where the orthogonal complement is taken in $\R^{n+1}$.
Thus, $\widetilde{V}_k^{\perp} =\text{span}\{\widetilde{\bm{e}}_k,\cdots,\widetilde{\bm{e}}_{n-1},\bm{u}\}$.
Note that $\dim \widetilde{V}_k = k$.

We first prove
\begin{equation}\label{inequality1}
\widetilde{\lambda}_k\leq \lambda_k + O(t),\quad 1\leq k\leq n-1.
\end{equation}
In fact, the min-max theorem brings about
\begin{eqnarray*}
\widetilde{\lambda}_k
&=& \min_{\substack{\widetilde{U} \subset \R^{n+1}\\ \dim \widetilde{U} = k+1}} \max_{\widetilde{\bm{v}} \in \widetilde{U} \setminus \{ 0 \}}
\frac{\langle\widetilde{A}\widetilde{\bm{v}},\widetilde{\bm{v}}\rangle}{\langle\widetilde{\bm{v}},\widetilde{\bm{v}}\rangle}
\,\,\leq\,\, \max_{\widetilde{\bm{v}}\in\widetilde{V}_{k+1}\setminus\{ 0\}}
\frac{\langle\widetilde{A}\widetilde{\bm{v}},\widetilde{\bm{v}}\rangle}{\langle\widetilde{\bm{v}},\widetilde{\bm{v}}\rangle}\\
&=& \max_{\bm{v}\in V_{k+1}\setminus\{ 0\}}
\frac{\langle A\bm{v},\bm{v}\rangle}{\langle \bm{v},\bm{v}\rangle} + O(t)
\,\,=\,\, 
\lambda_k + O(t) .
\end{eqnarray*}

Next we prove the reverse inequality $\widetilde{\lambda}_k\geq \lambda_k + O(t)$.
For any $\bm{v}\in V_k^{\perp}$, the squared Euclidean norm of
$\widetilde{\bm{w}} =\widetilde{\bm{v}} +c\bm{u} \in \widetilde{V}_k^{\perp}$, where $c$ is a constant, is
\begin{equation}\label{norm-compu}
\langle\widetilde{\bm{w}},\widetilde{\bm{w}}\rangle^2
= \langle\widetilde{\bm{v}},\widetilde{\bm{v}}\rangle^2 +c^2\langle\bm{u},\bm{u}\rangle^2\\
= \langle \bm{v},\bm{v}\rangle^2 + c^2 + \frac{v_n^2+c^2}{\alpha} t.
\end{equation}
By the max-min theorem, we have
$$
\widetilde{\lambda}_k
= \max_{\substack{\widetilde{W} \subset \R^{n+1}\\ \dim \widetilde{W} = k}} \min_{\widetilde{\bm{w}}\in \widetilde{W}^\perp}
\frac{\langle\widetilde{A}\widetilde{\bm{w}},\widetilde{\bm{w}}\rangle}{\langle\widetilde{\bm{w}},\widetilde{\bm{w}}\rangle}
\geq\min_{\widetilde{\bm{w}}\in\widetilde{V}_k^{\perp}}
\frac{\langle\widetilde{A}\widetilde{\bm{w}},\widetilde{\bm{w}}\rangle}{\langle\widetilde{\bm{w}},\widetilde{\bm{w}}\rangle}.
$$
For each small $t>0$, choose $\widetilde{\bm{w}}$, $\langle \widetilde{\bm{w}}, \widetilde{\bm{w}}\rangle=1$, which attains
the minimum on the rightmost side, and write $\widetilde{\bm{w}} =\widetilde{\bm{v}} +c(t) \bm{u}$.
By \eqref{norm-compu}, we have $|c(t)|\leq 1$ and $\langle \bm{v},\bm{v} \rangle \leq 
1$.
Writing $c:=c(t)$, we unravel
\begin{eqnarray}\label{lower-estimate-of-widetilde-lambda-k}
\widetilde{\lambda}_k
&\geq&\langle\widetilde{A}\widetilde{\bm{w}},\widetilde{\bm{w}}\rangle
\,\,=\,\, \langle\widetilde{A}\widetilde{\bm{v}},\widetilde{\bm{v}}\rangle +2c\langle\widetilde{A} \widetilde{\bm{v}}, \bm{u} \rangle + c^2\langle\widetilde{A}\bm{u},\bm{u}\rangle\nonumber\\
&=& \left( \langle A\bm{v},\bm{v}\rangle + O(t) \right)
+ \left( \frac{cv_n}{\alpha^{3/2}t^{1/2}} + O(t^{1/2}) \right)
+ \left( \frac{c^2}{t^2} + \frac{2c^2}{\alpha t} + O(t) \right)\\
&=& \langle A\bm{v},\bm{v}\rangle + \frac{c^2}{t^2} + \frac{2c^2}{\alpha t}
+ \frac{cv_n}{\alpha^{3/2}t^{1/2}} + O(t^{1/2}).\nonumber
\end{eqnarray}
The rightmost side is bounded below by $c^2/t^2+ O(1)$, since
$$
\frac{2c^2}{\alpha t} + \frac{cv_n}{\alpha^{3/2}t^{1/2}}
\geq \frac{2c^2}{\alpha t} -\frac{1}{2} \left( \frac{c^2}{\alpha t} + \frac{v_n^2}{\alpha^2} \right)
\geq - \frac{v_n^2}{2\alpha^2}.
$$
Combining this with \eqref{inequality1}, we obtain $c(t)=O(t)$,
and $\langle\bm{v},\bm{v}\rangle = 1 + O(t)$ by \eqref{norm-compu}.
Therefore, by \eqref{lower-estimate-of-widetilde-lambda-k},
\begin{equation}\label{inequality2}
\widetilde{\lambda}_k \geq \lambda_k + \frac{c^2}{t^2} + O(t^{1/2}),
\end{equation}
which clearly implies $\widetilde{\lambda}_k \geq \lambda_k + O(t^{1/2})$.
This concludes the proof.
\end{proof}

\begin{Remark}
The dominant part of the $2 \times 2$ block on the bottom right corner of the matrix
$\widetilde{A}$ is
$$
S=
\begin{pmatrix}
\frac{t^{-1}}{\alpha} & -\frac{t^{-3/2}}{\sqrt{\alpha}}\\
-\frac{t^{-3/2}}{\sqrt{\alpha}} & t^{-2}\\
\end{pmatrix},
$$
and $S$ is most crucial for the divergence of the largest eigenvalue of $\widetilde{A}$.
The matrix $S$ has the eigenvalues $0$ and $(\alpha+t)/(\alpha t^2)$ with corresponding
eigenvectors ${}^t (1, \sqrt{t/\alpha})$ and ${}^t (-\sqrt{t/\alpha}, 1)$, respectively.
In the proof of Proposition \ref{asymptotic}, the definition of $\widetilde{\bm{v}}$
and the choice of $\bm{u}$ are motivated by these eigenvectors of $S$.
\end{Remark}

\begin{Remark}
It follows from \eqref{inequality1}, \eqref{inequality2} that $c(t)$ is at most $O(t^{5/4})$.
\end{Remark}

Theorem \ref{contract} follows by applying Proposition \ref{asymptotic} to the Laplacian
matrices of $(G,l)$ and $(\widetilde{G},l_t)$.

\section{Monotonicity}
In this section, we prove that the invariant $\Lambda_1(G)$ is nonincreasing under taking a subgraph. 
Namely, we prove 
\begin{Theorem}
If $G$ is a finite connected graph and $H$ is a subgraph of $G$, then the inequality 
$\Lambda_1(H)\leq \Lambda_1(G)$ holds.
\end{Theorem}

\begin{proof}
Choose a vertex $u$ such that the distance between $u$ and $H$ is maximum.
If $\mathrm{deg}(u)=1$, by contracting the edge $e$ with end point $u$, we obtain
a graph $G'$, which is connected and satisfies $\Lambda_1(G')\leq \Lambda_1(G)$
by Theorem \ref{contract}.
If $d:=\mathrm{deg}(u)\geq 2$, let $e_1,\dots, e_d$ denote the edges with end point $u$
and cut the graph $G$ at $u$ to meke $e_1$ a pendant edge, that is, an edge incident to a degree-one vertex. 
Otherwise said, we introduce a clone vertex $\overline{u}$ and re-join $e_1$ to $\overline{u}$ instead of $u$.
We denote the graph thus obtained by $G_{\mathrm{cut}}$. Note $G_{\mathrm{cut}}$ is connected and we identify 
$E(G_{\mathrm{cut}})\cong E(G)$.
Then for any length function $l\colon E(G)\to \R_{>0}$, we have

\medskip\noindent
{\bf Claim.}
\begin{equation}\label{duplicate}
\lambda_1(G_{\mathrm{cut}}, l)\leq \lambda_1(G, l).
\end{equation}

\medskip\noindent
{\it Proof of Claim.}\quad
Recall that
$$
\lambda_1(G,l) = \inf_{\varphi\perp {\bf 1}} \frac{\sum_{vw\in E(G)} m_1(vw)
(\varphi(v)-\varphi(w))^2}{\sum_{v\in V(G)} m_0(v)\varphi(v)^2}.
$$
Extend any $\varphi\colon V(G)\to \R$ such that $\varphi\perp {\bf 1}$ in $\R^{V(G)}$
to $\widetilde{\varphi}\colon V(G_{\mathrm{cut}})\to \R$ by setting $\widetilde{\varphi}(\overline{u})=\varphi(u)$.
Then it is easy to see that $\widetilde{\varphi}\perp {\bf 1}$ in $\R^{V(G_{\mathrm{cut}})}$ and
the numerator/denominator of the Rayleigh quotient for $\widetilde{\varphi}$ with respect to
$(G_{\mathrm{cut}}, l)$ are the same as those above.
Since the collection of all functions $\psi$ such that $\psi\perp {\bf 1}$
in $\R^{V(G_{\mathrm{cut}})}$ is larger that that of all $\widetilde{\varphi}$ given as above, we conclude
the desired result.

\medskip
As an immediate consequence of the claim, we have $\Lambda_1(G_{\mathrm{cut}})\leq \Lambda_1(G)$.
Then by contracting the edge $e_1$ in $G_{\mathrm{cut}}$, we obtain a graph ${G_{\mathrm{cut}}^{\mbox{}}}'$, which is connected
and satisfies $\Lambda_1({G_{\mathrm{cut}}^{\mbox{}}}')\leq \Lambda_1(G_{\mathrm{cut}})$ by Theorem \ref{contract} again.
Therefore, $\Lambda_1({G_{\mathrm{cut}}^{\mbox{}}}')\leq \Lambda_1(G)$.

We repeat the above procedure until only the vertices of $H$ remain. 
We may still have an edge of $G$ which is not an edge of $H$. Such an edge is not pendant, but the above cutting argument still works. 
By repeating this latter procesure until only the edges of $H$ remain, the theorem is proved.
\end{proof}

\end{document}